\documentstyle[11pt]{article}
\input amssym.def
\input amssym.tex
\input psfig
\makeatletter
\def\@sect#1#2#3#4#5#6[#7]#8{\ifnum #2>\c@secnumdepth
     \def\@svsec{}\else 
     \refstepcounter{#1}\edef\@svsec{\csname the#1\endcsname.\hskip .75em }\fi
     \@tempskipa #5\relax
      \ifdim \@tempskipa>\z@ 
        \begingroup #6\relax
          \@hangfrom{\hskip #3\relax\@svsec}{\interlinepenalty \@M #8\par}%
        \endgroup
       \csname #1mark\endcsname{#7}\addcontentsline
         {toc}{#1}{\ifnum #2>\c@secnumdepth \else
                      \protect\numberline{\csname the#1\endcsname}\fi
                    #7}\else
        \def\@svsechd{#6\hskip #3\@svsec #8\csname #1mark\endcsname
                      {#7}\addcontentsline
                           {toc}{#1}{\ifnum #2>\c@secnumdepth \else
                             \protect\numberline{\csname the#1\endcsname}\fi
                       #7}}\fi
     \@xsect{#5}}
\def\@begintheorem#1#2{\it \trivlist \item[\hskip \labelsep{\bf #1\ #2.}]}
\makeatother

\newtheorem{theorem}{Theorem}
\newtheorem{lemma}{Lemma}

\def\slashedfrac#1#2{\hbox{\kern.1em %
 \raise.5ex\hbox{\the\scriptfont0 #1}\kern-.11em %
 /\kern-.15em\lower.25ex\hbox{\the\scriptfont0 #2}}}

\newcommand{\hsp}{\hspace*{\parindent}}

\newcommand{\eeq}{\end{equation}}
\newcommand{\beql}[1]{\begin{equation}\label{#1}}
\newcommand{\bsq}{{\vrule height .9ex width .8ex depth -.1ex }}

\newcommand{\RR}{{\Bbb R}}

\newcommand{\sA}{{\cal A}}
\newcommand{\La}{\Lambda}

\newcommand{\dd}{\ldots}

\begin{document}
\begin{center}
{\Large {\bf The Antipode Construction for Sphere Packings}} \\
\vspace{\baselineskip}
{\em J.~H. Conway} \\
\vspace{.25\baselineskip}
Mathematics Department \\
Princeton University \\
Princeton, New Jersey 08544 \\
\vspace{1\baselineskip}
{\em N.~J.~A. Sloane} \\
\vspace{.25\baselineskip}
Mathematical Sciences Research Center \\
AT\&T Bell Laboratories \\
Murray Hill, New Jersey 07974 \\
\vspace{1.5\baselineskip}
{\bf Summary}
\vspace{.5\baselineskip}
\end{center}
\setlength{\baselineskip}{1.5\baselineskip}

A construction for sphere packings is introduced that is
parallel to the ``anticode''
construction for codes.
This provides a simple way to view Vardy's recent 20-dimensional sphere packing,
and also produces packings in dimensions 22, 44--47 that are denser than any previously known.

\clearpage
\large\normalsize
\renewcommand{\baselinestretch}{1}
\thispagestyle{empty}
\setcounter{page}{1}
\begin{center}
{\Large {\bf The Antipode Construction for Sphere Packings}} \\
\vspace{\baselineskip}
{\em J.~H. Conway} \\
\vspace{.25\baselineskip}
Mathematics Department \\
Princeton University \\
Princeton, New Jersey 08544 \\
\vspace{1\baselineskip}
{\em N.~J.~A. Sloane} \\
\vspace{.25\baselineskip}
Mathematical Sciences Research Center \\
AT\&T Bell Laboratories \\
Murray Hill, New Jersey 07974 \\
\vspace{1.5\baselineskip}
\end{center}
\setlength{\baselineskip}{1.5\baselineskip}
\section{Introduction}
\hsp
Vardy \cite{Vardy} recently found a remarkable
20-dimensional sphere packing that is denser than
any previously known, and asked if the affine automorphism group of the packing
is transitive on spheres.
We will give an alternative simple construction which yields
Vardy's packing and some other new packings.
We will also show that the group of Vardy's packing does act
transitively.

For undefined terms from sphere packing, see \cite[Chap.~1]{CS93}.
\section{The antipode construction}
\hsp
The following construction is an analogue of the ``anticode'' construction
for codes \cite[Chap.~17, \S6]{MS77}.
(The common theme of these two constructions is that instead of looking for well-separated
points, which is what most constructions do,
now we look for points somewhere else that are close together and factor them out.)

Let $\Lambda$ be a lattice with minimal norm (i.e. squared
length) $\mu$ in an $n$-dimensional Euclidean space $W$.
Let $U,V$ be respectively $k$- and $l$-dimensional vector spaces
with $W = U \oplus V$, $n=k+l$, such that
$\Lambda \cap U=K$ and $\Lambda \cap V=L$ are $k$- and $l$-dimensional lattices,
respectively.
The projection maps from $W$ to $U$ and $V$ will be denoted by $\pi_U$, $\pi_V$, respectively,
and we define the lattices
$M := \pi_U ( \Lambda )$,
$N := \pi_V ( \Lambda )$.
In the examples below, $\Lambda$ will be self-dual, and so
$M$ and $N$
will be the lattices $K^\ast$
and $L^\ast$ dual to $K$ and $L$
\cite[p.~166]{CS93}.

Suppose we can find a subset $S = \{ u_1 , \dd , u_s \}$ of $M$ such that
$dist^2 (u_i , u_j) \le \beta$ for all $i,j$,
for some constant $\beta$.
We define $\sA (S)$ to consist of the points
$$\{ \pi_V (w) : ~ w \in \Lambda ,~ \pi_U (w) \in S \} ~.$$
We say that the sets $S$ and $\sA (S)$ are
{\em antipodal} to each other - this explains our name for this construction.

\begin{theorem}
\label{th1}
By drawing spheres of radius $\rho = \frac{1}{2} \sqrt{\mu - \beta}$ around the points of $\sA (S)$ we obtain an $l$-dimensional sphere packing of center density
\beql{eq1}
\delta = \frac{s \cdot \rho^l}{\sqrt{\det L}} ~.
\eeq
\end{theorem}

\noindent{\bf Proof.}
$\sA (S)$ is the union of $s$ translates of $L$, and by construction the squared
distance between the points of $\sA (S)$ is at least $\mu - \beta$.
The result follows from Eq.~(13) of \cite[Chap.~1]{CS93}.~~~$\bsq$

In many cases (in particular, in dimensions just below 8, 24 and 48) the densest
packings known up to now are obtained by taking a self-dual
lattice $\Lambda$ $(E_8$, $\Lambda_{24}$ and
$P_{48p}$ or $P_{48q}$, respectively) and using suitable
sections $L= \Lambda \cap V$.
Some of these densities can now be improved by the antipode
construction, since
$\sA (S)$ is denser than $L$ whenever
\beql{eq2}
s \left( 1- \frac{\beta}{\mu} \right)^{\frac{l}{2}} > 1 ~,
\eeq
the left-hand side of this formula
being the ratio of their densities.

This construction improves on the old records in the following cases.
In the first two examples we take $n=24$ and $\Lambda$ to be the Leech lattice
$\Lambda_{24}$ in $\RR^{24}$, with minimal norm $\mu =4$.
We write the vectors of $\Lambda_{24}$ in the form
$\frac{1}{\sqrt{8}} (a_1 , \dd , a_{24})$, where the $a_i$ are the coordinates of a MOG diagram
\cite[Chap.~11]{CS93}.
\paragraph{Vardy's 20-dimensional packing \protect\cite{Vardy}.}
We take $l= 20$, $V =$~vectors in $\RR^{24}$
beginning with four zeros, $L =$~laminated lattice $\Lambda_{20}$ with determinant 64,
$K = \sqrt{2} D_4$ (the root lattice $D_4$, rescaled so as to have minimal norm 4),
$M= \frac{1}{\sqrt{2}} D_4^\ast$ (the dual lattice $D_4^\ast$, rescaled so as to have
minimal norm $\frac{1}{2}$),
$s=4$, $S= \frac{1}{\sqrt{8}} \{0~0~0~0,~~1~1~1~1,~~2~0~0~0,~~1~1~1~-1 \}$,
$\beta = \frac{1}{2}$, obtaining a packing $V_{20}$ with center density
$$\frac{7^{10}}{2^{31}} = .1315 \dd ~,$$
compared with the old record of $1/8 = .125$ for $\Lambda_{20}$.
It is easy to show that this construction is equivalent to Vardy's (compare Eq.~(9) of \cite{Vardy}).
Vardy shows that each sphere touches 15360 others.
\paragraph{Dimension 22.}
Again using the Leech lattice, we take $l=22$, $V =$~vectors in $\RR^{24}$ beginning
with three equal coordinates, $L = \Lambda_{22}$ with determinant 12,
$K= \sqrt{2} A_2$, $M= \frac{1}{\sqrt{2}} A_2^\ast$,
and $S$ to consist of the $s=3$ vertices of an equilateral triangle in $M$, with
$\beta = \frac{1}{3}$, obtaining a packing $V_{22}$ with center density
$$\delta = \frac{11^{11}}{2^{23} 3^{10.5}} = .3325 \dd ~,$$
a new record, compared with $1/ \sqrt{12} = .2886 \dd$ for $\Lambda_{22}$.

Explicit coordinates for $V_{22}$ are as follows.
The three translates of $\La_{22}$ may be obtained by taking all vectors of $\La_{24}$ in which the first three coordinates have the form $\frac{1}{\sqrt{8}} (a,a,a)$,
\linebreak
$\frac{1}{\sqrt{8}} (a+2, a,a)$, or
$\frac{1}{\sqrt{8}} (a,a-2,a)$.
Their projections onto $U$ are then proportional to
$(0,0,0)$,
$\left( \frac{2}{3}, \,- \frac{1}{3}, \,- \frac{1}{3} \right)$,
and $\left( \frac{1}{3}, \,- \frac{2}{3}, \,\frac{1}{3} \right)$,
which do indeed form an equilateral triangle.
The corresponding vectors in $V_{22}$, which are the projections onto $V$, are
$\frac{1}{\sqrt{8}} (a,a,a, \dd )$, $\frac{1}{\sqrt{8}} \left( a + \frac{2}{3} , a + \frac{2}{3} , a+ \frac{2}{3} , \dd \right)$, and
$\frac{1}{\sqrt{8}} \left( a- \frac{2}{3} , a - \frac{2}{3} , a- \frac{2}{3} , \dd \right)$.

As we will see in \S3, the affine groups of both $V_{20}$ and $V_{22}$ act transitively on the spheres.
To find the kissing number of $V_{22}$, consider the sphere centered at the zero vector.
No sphere in this copy of $\La_{22}$ touches it.
The number in the second translate that touch it is equal to the number of minimal vectors
in $\La_{24}$ which begin
$\frac{1}{\sqrt{8}} (a+2,a,a, \dd )$.
These vectors are
$$\begin{array}{lll}
\mbox{for}~~a=-2: & 0, -2, -2; \pm 2^6 0^{15} & (\# = 2^5 \cdot 56) ~, \\ [+.1in]
\mbox{for}~~a= -1: & 1, -1, -1; \pm 1^{20} \pm 3 & (\# = 2^9 \cdot 21 ) ~, \\ [+.1in]
\mbox{for}~~a=0 : & 2,0,0 ; \pm 2^7 0^{14} & (\# = 2^6 \cdot 120 ) ~, \\ [+.1in]
\mbox{for}~~a=1: & 3,1,1; \pm 1^{21} & (\# = 2^9 \cdot 1 )
\end{array}
$$
(cf. \cite{CS93}, p.~278), for a total of 20736.
Therefore each sphere in this packing touches 41472 others.
\paragraph{Dimensions 44--47.}
We take $n=48$, and $\Lambda =$ either $P_{48p}$ or $P_{48q}$ \cite[p.~195]{CS93}
(so that there are at least two
inequivalent packings in each of these cases), with $\mu =6$.
As on p.~168 of \cite{CS93} we can find subspaces $U$ in $\RR^{48}$ of dimensions 1, 2, 3, 4 such that $M$ is
respectively $\frac{1}{\sqrt{3}} A_1^\ast$, $\frac{1}{\sqrt{3}} A_2^\ast$,
$\frac{1}{\sqrt{3}} A_3^\ast \cong \frac{1}{\sqrt{3}} D_3^\ast$, $\frac{1}{\sqrt{3}} D_4^\ast$.
In these four lattices it is easy to find $s=2,3,4,4$ points, respectively,
for which $\beta = \frac{1}{6}$, $\frac{2}{9}$, $\frac{1}{3}$, $\frac{1}{3}$.
We then obtain
packings in dimensions $l= 47, 46, 45, 44$ with center
densities $2^{-70} 3^{-24} 35^{23.5} = 5788.8 \dd$,
$3^{-46.5} 13^{23} = 2719.9 \dd$, $2^{-44} 3^{-24} 17^{22.5} = 974.6 \dd$,
$2^{-43} 3^{-24} 17^{22} = 472.7 \dd$, which
improve on the old records by factors of $1.031 \dd$, $1.259 \dd$, $1.105 \dd$ and 1.137,
respectively.

We found no other improvements on the present records for density
(cf. \cite{CS93}, Chap.~1), despite examining all likely cases in
low dimensions.
\section{The group of Vardy's packing}
\hsp
We will determine $\mbox{Aff}(V_{20})$, the full affine symmetry group of $V_{20}$ and in particular
show that it acts transitively on the spheres, thus answering a question raised in \cite{Vardy}.
Similar but easier arguments apply to $V_{22}$.
\begin{lemma}
\label{le1}
Any affine symmetry of $V_{20}$ is induced by an affine symmetry of $\La_{24}$.
\end{lemma}

\noindent{\bf Proof.}
$V_{20}$ consists of four translates of $\La_{20}$, which we shall label $[0]$, $[1]$, $[2]$,
$[3]$, corresponding to the four members of $S$.
Two points of $V_{20}$ have integral squared distance if and only if they are in the same
translate, so this decomposition into translates is intrinsic.
In an abstract 4-space we take four vectors
$v_0 = 0~0~0~0$, $v_1 = 1~1~1~1$, $v_2 = 2~0~0~0$, $v_3 =1~1~1~-1$,
and glue $v_i$ to $[i]$ (cf. \cite[Chap.~4]{CS93}).
The lattice affinely generated by the result is the laminated lattice $\La_{23}$,
consisting of those vectors of $\La_{24}$ with second and third coordinates equal.
So any affine symmetry of $V_{20}$ extends uniquely to an affine symmetry of $\La_{23}$.
A similar argument shows that any affine symmetry of $\La_{23}$ extends uniquely to an affine symmetry of $\La_{24}$.
This completes the proof.~~~$\bsq$

So we may regard $\mbox{Aff}(V_{20})$ as a subgroup of the affine symmetry group $\mbox{Co}_\infty$ of $\La_{24}$.
In fact $\mbox{Aff} (V_{20})$ is the subgroup of $\mbox{Co}_\infty$ that preserves the partition of $W= \RR^{24}$
into $U \oplus V$, and also fixes (setwise) the vertices of the tetrahedron formed by the elements
of $S$.
It is easy to see that this group acts as a symmetric group ${\bf S}_4$ on the vertices of the tetrahedron.
\begin{lemma}
\label{le2}
$\mbox{Aff}(V_{20})$ acts transitively on $V_{20}$.
\end{lemma}

\noindent{\bf Proof.}
By the previous remark, $\mbox{Aff}(V_{20})$ acts as an ${\bf S}_4$ on the four cosets of $\La_{20}$.
But $\mbox{Aff} (V_{20})$ also includes translates by vectors of $\La_{20}$,
which act transitively in each coset.~~~$\bsq$
\begin{lemma}
\label{le3}
The subgroup of $\mbox{Aff} (V_{20})$ that fixes one sphere is a group of order
$2^8 \, | {\bf M}_{20} | \, 2 \cdot 3! = 2949120$,
and acts transitively on the 15360 neighboring spheres.
\end{lemma}

We omit the proof, which is a straightforward calculation using the results of
\cite[Chap.~11]{CS93}.
The factor $2^8 | {\bf M}_{20} |$ is the order of the subgroup of $2^{12} {\bf M}_{24}$ that fixes the first four coordinates,
the factor 2 corresponds to the interchange of the second and third
coordinates,
and the factor 3! refers to permutations of the remaining
three vertices of the tetrahedron.

In other words, not only does every sphere in Vardy's
packing look like every other one, but every pair of
spheres in contact looks like every other such pair.
However, there is more than one type of triple of mutually contiguous spheres.
\paragraph{Acknowledgments.}
We are grateful to Alexander Vardy for sending us a preprint of \cite{Vardy}.

\end{document}